\documentclass[11pt,onecolumn]{article}
\usepackage[a4paper,left=1.7cm,right=1.7cm,top=2.2cm,bottom=2.5cm,columnsep=0.8cm]{geometry}
\usepackage{graphicx}
\usepackage{amsmath,amsfonts}
\usepackage{amssymb}
\usepackage{algorithm}
\usepackage{algpseudocode}
\usepackage{booktabs}
\usepackage{url}
\usepackage{natbib}
\usepackage[subrefformat=parens]{subcaption}
\newenvironment{proof}[1][Proof]{\textbf{#1.} }{\hfill$\square$}
\newtheorem{theorem}{Theorem}
\newtheorem{lemma}{Lemma}
\newtheorem{remark}{Remark}
\newtheorem{prob}{Problem}
\newtheorem{ass}{Assumption}

\newcommand{\A}{{\mathcal A}}
\newcommand{\B}{{\mathcal B}}

\newcommand{\C}{{\mathcal C}}

\title{Optimal Dynamic Quantizers for Periodically Time-Varying Systems via Cyclic Reformulation}

\author{Hiroshi Okajima\thanks{Faculty of Advanced Science and Technology, Kumamoto University, Japan. E-mail: \texttt{okajima@mech.kumamoto-u.ac.jp}}}

\date{}

\begin{document}

\maketitle

\begin{abstract}
This paper addresses the design of optimal dynamic quantizers for discrete-time linear periodically time-varying (LPTV) plants with discrete-valued inputs. Optimal dynamic quantizers, which minimize the worst-case output error caused by quantization, are well established for linear time-invariant (LTI) plants, whereas the corresponding design for plants with periodically time-varying dynamics has not been addressed. For square LPTV plants having relative degree one at every time phase, this paper derives an optimal periodic dynamic quantizer in closed form. The quantizer is constructed by applying the LTI theory to the cycled plant obtained by cyclic reformulation and reading off its periodic coefficient matrices block by block; its optimality is then verified directly, showing that it attains a performance lower bound determined by the largest instantaneous input-to-output gain over the time phases. The worst-case output error of general periodic dynamic quantizers is also characterized. The optimal quantizer is shown to be stable if and only if all invariant zeros of the cycled plant lie in the open unit disk, which is the condition that governs stable causal inversion of LPTV plants. A numerical example with comparisons to conventional quantizers illustrates the results.
\end{abstract}


\section{Introduction}\label{sec:intro}

Control systems with discrete-valued input signals arise in many practical situations, including mechanical systems with ON/OFF actuators, networked control systems with finite communication bandwidth, and power electronics with switching converters. The design of quantizers that map continuous-valued control signals to discrete-valued ones is a central issue in such systems.

For linear time-invariant (LTI) plants, \citet{azuma2008automatica,azuma2008tac} established a systematic framework for the design of optimal dynamic quantizers, which minimize the worst-case output error caused by quantization. In particular, \citet{azuma2008automatica} derived a closed-form optimal dynamic quantizer that cancels the influence of past quantization errors on the output, leaving only the instantaneous error determined by the input-to-output gain $CB$. The stability of this quantizer is governed by the zeros of the plant; for nonminimum phase plants, \citet{azuma2008tac} developed a linear programming-based synthesis, and \citet{sawada2011acc} proposed an LMI-based invariant set approach for stable suboptimal quantizer design.

The framework has also been extended to nonlinear plants \citep{azuma2012tac} and to decentralized quantizer structures \citep{minami2009acc}, and \citet{fujimoto2026ifacsc} recently proposed a data-driven design that avoids the need for an explicit plant model. Dynamic quantizers have been applied to feedback control \citep{minami2007cdc,azuma2008cdc}, to networked control systems with communication rate constraints \citep{okajima2016tac,almakhles2017tcst} or with time-varying delays and packet losses \citep{li2015ietcta}, and to encrypted control \citep{teranishi2021sice}.

Quantization and discrete-valued control have also been studied in other contexts. From a signal processing perspective, \citet{quevedo2005tcsi} proposed a multistep optimal converter that generalizes $\Delta\Sigma$ modulation, and \citet{ohno2017tcsi} addressed noise shaping filter optimization for quantizers. In quantized feedback stabilization, quantizers whose sensitivity is adjusted online were introduced by \citet{brockett2000tac}, and quantization with dynamic or encoding--decoding mechanisms continues to be studied in iterative learning control \citep{liu2025tcyb} and in state estimation over networks \citep{bai2025jiot,meng2023tnse}. Finally, discrete-valued control has been approached through sum-of-absolute-values optimization \citep{ikeda2017tac}.

Despite these advances, optimal dynamic quantizer design for plants with periodically time-varying dynamics has not been addressed. Linear periodically time-varying (LPTV) systems arise in multirate sampled-data systems \citep{chen1995sampled}, in systems with periodic scheduling of resources, and in plants subject to periodic operating conditions \citep{vanzundert2019ijc}. When the actuator of such a system accepts only discrete-valued inputs, the quantizer must be designed for a periodic plant. In a preliminary study, \citet{okajima2019sice} proposed a lifting-based periodically time-varying dynamic quantizer for networked control systems with periodically varying communication rate constraints, enabling effective quantization at non-integer bit rates such as 1.5 bits per sampling period. However, that work addressed a time-varying quantizer structure under fixed LTI plant dynamics. Dynamic quantizers for multirate control systems were also studied by \citet{morita2011iscie}, but for an LTI plant combined with hold devices of different periods.

This paper derives a closed-form optimal periodic dynamic quantizer for square LPTV plants having relative degree one at every time phase, i.e., plants whose instantaneous input-to-output gain is nonsingular at every phase; this is the periodic counterpart of the LTI result of \citet{azuma2008automatica}. The quantizer is obtained through cyclic reformulation \citep{cyc1,cyc2}, which transforms an LPTV system into an equivalent LTI system at the original sampling rate while preserving the step-by-step input--output correspondence. This technique has been applied to robust estimator synthesis \citep{hosoe2019sice}, system identification \citep{okajima1}, and stable inversion of LPTV plants \citep{okajima2026inv}. Applying the LTI theory to the cycled plant yields an LTI quantizer whose coefficient matrices retain the cyclic block structure, so that a periodically time-varying quantizer with the original state dimension can be read off block by block; its optimality is then verified directly in the original time domain.

The contributions of this paper are as follows. First, the optimal dynamic quantizer design problem for LPTV plants is formulated, and a closed-form optimal periodic quantizer is derived together with the minimum achievable worst-case output error, generalizing the LTI result of \citet{azuma2008automatica}. Second, the worst-case output error of a general periodic dynamic quantizer is characterized in terms of phase-wise impulse response matrices. Third, the optimal quantizer is shown to be stable if and only if the plant is periodically minimum phase, i.e., all invariant zeros of the cycled plant lie in the open unit disk, which is exactly the condition for stable causal inversion of the plant \citep{okajima2026inv}.

This paper is organized as follows. Section~\ref{sec:methods} presents the preliminaries on cyclic reformulation, reviews the existing optimal dynamic quantizer theory for LTI systems, and states the problem formulation for LPTV systems. Section~\ref{sec:results} presents the main results: the performance analysis of general periodic quantizers, the optimal quantizer, its stability analysis, and a numerical example. Section~\ref{sec:discussion} discusses the role of cyclic reformulation and the limitations of the results, and Section~\ref{sec:conclusion} concludes the paper.

\textit{Notation:}
Let $\mathbb{R}$, $\mathbb{R}_+$, $\mathbb{Z}$, $\mathbb{Z}_+$, and $\mathbb{N}$ denote the real numbers, positive reals, integers, nonnegative integers, and positive integers, respectively. For a matrix $M = \{M_{ij}\}$, $\mathrm{abs}(M):=\{|M_{ij}|\}$ denotes the elementwise absolute value, and $\|M\|$ denotes the $l_\infty$-induced norm (maximum absolute row sum). For a square matrix $M$, $\rho(M)$ denotes the spectral radius, and $M$ is called Schur stable if $\rho(M) < 1$. For a vector $x$, $\|x\|$ is the $l_\infty$-norm. For a signal sequence $X = \{x_0, x_1, \ldots\}$, $\|X\|:=\sup_{k} \|x_k\|$; the space $l_\infty^m$ denotes the set of bounded sequences in $\mathbb{R}^m$. The symbols $I_n$ and $0_{n\times m}$ denote the $n \times n$ identity matrix and the $n \times m$ zero matrix (the subscripts are omitted when the dimensions are clear from the context), $\mathbf{1}_n$ is the $n$-vector of ones, and $\otimes$ denotes the Kronecker product. For a square system $(A, B, C)$ with zero feedthrough, the invariant zeros are the values $z_0 \in \mathbb{C}$ at which the rank of the system pencil $\left[\begin{smallmatrix} z_0 I - A & -B \\ C & 0 \end{smallmatrix}\right]$ drops below its normal rank. The static quantizer $q:\mathbb{R}^m \to \mathbb{V}^m$ with $\mathbb{V}:=d\mathbb{Z}$ and quantization interval $d \in \mathbb{R}_+$ is the nearest-neighbor quantizer satisfying $\mathrm{abs}(q(\mu)-\mu) \leq \frac{d}{2}\mathbf{1}_m$ for all $\mu \in \mathbb{R}^m$. Throughout this paper, all subscripts of $N$-periodic matrices (such as $A_k$, $B_k$, $C_k$, and their derived quantities) are interpreted modulo $N$ without further notice; that is, $A_k = A_{k \bmod N}$ (e.g., $A_{-1} = A_{N-1}$ and $C_{N+1} = C_1$).

\section{Preliminaries and Problem Formulation}\label{sec:methods}

\subsection{Cyclic Reformulation of Periodically Time-Varying Systems}\label{sec:cyclic}

Consider a discrete-time LPTV system with period $N$ ($N \in \mathbb{N}$):
\begin{eqnarray}
x(k+1) &=& A_k x(k) + B_k u(k), \label{eq:lptv_state}\\
y(k) &=& C_k x(k), \label{eq:lptv_output}
\end{eqnarray}
where $x \in \mathbb{R}^n$, $u \in \mathbb{R}^m$, $y \in \mathbb{R}^p$, and $A_{k+N} = A_k$, $B_{k+N} = B_k$, $C_{k+N} = C_k$ for all $k \in \mathbb{Z}_+$. Here, $D_k = 0$ is assumed for simplicity, which is natural for many sampled-data systems.

The cyclic reformulation \citep{cyc1,cyc2} transforms this LPTV system into a time-invariant system. Letting $i := k \bmod N \in \{0, 1, \ldots, N-1\}$ denote the time phase, the cycled input signal $\check{u}(k) \in \mathbb{R}^{Nm}$ is constructed by placing $u(k)$ in the $i$-th block position and setting all other blocks to zero:
\begin{equation}\label{eq:cycled_input}
\check{u}(k) = e_{i} \otimes u(k), \quad i := k \bmod N,
\end{equation}
where $e_i \in \mathbb{R}^N$ is the $(i+1)$-th standard unit vector. For example, $\check{u}(0) = [u(0)^\top, 0, \ldots, 0]^\top$, $\check{u}(1) = [0, u(1)^\top, 0, \ldots, 0]^\top$, and so on. The cycled state $\check{x}(k) \in \mathbb{R}^{Nn}$ and cycled output $\check{y}(k) \in \mathbb{R}^{Np}$ are defined analogously.

The cycled system is the following time-invariant system:
\begin{eqnarray}
\check{x}(k+1) &=& \check{A}\check{x}(k) + \check{B}\check{u}(k), \label{eq:cycled_state}\\
\check{y}(k) &=& \check{C}\check{x}(k), \label{eq:cycled_output}
\end{eqnarray}
where $\check{A}$ and $\check{B}$ have a block-circulant structure (the nonzero block of the $j$-th block-column appears in block-row $j+1$) and $\check{C}$ is block-diagonal:
\begin{equation}\label{eq:checkA}
\check{A} = \begin{bmatrix}
0 & 0 & \cdots & 0 & A_{N-1} \\
A_0 & 0 & \cdots & 0 & 0 \\
0 & A_1 & \ddots & \vdots & \vdots \\
\vdots & \ddots & \ddots & 0 & 0 \\
0 & \cdots & 0 & A_{N-2} & 0
\end{bmatrix} \in \mathbb{R}^{Nn \times Nn},
\end{equation}
\begin{equation}\label{eq:checkB}
\check{B} = \begin{bmatrix}
0 & 0 & \cdots & 0 & B_{N-1} \\
B_0 & 0 & \cdots & 0 & 0 \\
0 & B_1 & \ddots & \vdots & \vdots \\
\vdots & \ddots & \ddots & 0 & 0 \\
0 & \cdots & 0 & B_{N-2} & 0
\end{bmatrix} \in \mathbb{R}^{Nn \times Nm},
\end{equation}
\begin{equation}\label{eq:checkC}
\check{C} = \mathrm{diag}(C_0, C_1, \ldots, C_{N-1}) \in \mathbb{R}^{Np \times Nn}.
\end{equation}
With $\check{x}(0) = e_0 \otimes x(0)$ and the cycled input (\ref{eq:cycled_input}), the trajectories of (\ref{eq:cycled_state})--(\ref{eq:cycled_output}) satisfy $\check{x}(k) = e_i \otimes x(k)$ and $\check{y}(k) = e_i \otimes y(k)$ with $i = k \bmod N$ \citep[Chap.~6]{cyc2}. The cycled system thus reproduces the LPTV system step by step at the original sampling rate: at each time $k$, only the $i$-th block of each cycled signal is nonzero, and it equals the original signal.

\subsection{Optimal Dynamic Quantizers for LTI Systems}\label{sec:lti_review}

This subsection reviews the key results of \citet{azuma2008automatica,azuma2008tac} for LTI plants. Consider the system $\Sigma$ composed of an LTI plant $P$ and a dynamic quantizer $Q$. The plant is
\begin{equation}\label{eq:lti_plant}
P: \quad x(k+1) = Ax(k) + Bv(k), \quad y(k) = Cx(k),
\end{equation}
where $x \in \mathbb{R}^n$, $v \in \mathbb{R}^m$, $y \in \mathbb{R}^p$. The dynamic quantizer is of the form
\begin{equation}\label{eq:lti_quantizer}
Q: \begin{cases}
\xi(k+1) = \A\xi(k) + \B_1 u(k) + \B_2 v(k), \\
v(k) = q(\C\xi(k) + u(k)),
\end{cases}
\end{equation}
where $\xi \in \mathbb{R}^{n_\xi}$ with $\xi(0) = 0$ and $u$ is the continuous-valued input to $\Sigma$. The maximum output difference is defined as
\begin{equation}\label{eq:E_def}
E(Q) := \sup_{U \in l_\infty^m} \|Y(x_0, U) - Y(x_0, Q(U))\|,
\end{equation}
where $Y(x_0, U)$ and $Y(x_0, Q(U))$ denote the output sequences of $P$ with initial state $x_0$ driven by $U$ directly and through $Q$, respectively. The following lemma summarizes Theorems~1(i) and~2 of \citet{azuma2008automatica}.

\begin{lemma}\label{lem:lti_optimal}
For the system $\Sigma$, if $m = p$ and $CB$ is nonsingular, then
\begin{equation}
E(Q) \geq \|\mathrm{abs}(CB)\| \cdot \frac{d}{2}
\end{equation}
for any dynamic quantizer $Q$ of the form (\ref{eq:lti_quantizer}). Moreover, an optimal dynamic quantizer is given by
\begin{equation}\label{eq:Qc}
Q^c: \begin{cases}
\xi(k+1) = A\xi(k) + Bv(k) - Bu(k), \\
v(k) = q(-(CB)^{-1}CA\xi(k) + u(k)),
\end{cases}
\end{equation}
and the minimum value of $E(Q)$ is $E(Q^c) = \|\mathrm{abs}(CB)\| \cdot \frac{d}{2}$.
\end{lemma}
Since $\B_1 = -B$ and $\B_2 = B$ in $Q^c$ (i.e., $\B_1 = -\B_2$), the state equation reads $\xi(k+1) = A\xi(k) + B(v(k) - u(k))$, so that $\xi$ tracks the accumulated quantization effect in the plant state coordinates, and the term $-(CB)^{-1}CA\xi(k)$ compensates for the influence of past quantization errors on the next output $y(k+1)$.

The following result \citep[Theorem~1(ii)]{azuma2008automatica} characterizes the performance of a general dynamic quantizer.

\begin{lemma}\label{lem:lti_performance}
For the system $\Sigma$ composed of $P$ in (\ref{eq:lti_plant}) and $Q$ in (\ref{eq:lti_quantizer}) with $\B_1 = -\B_2$, let $w(k) := v(k) - (\C\xi(k) + u(k))$ denote the static quantization error satisfying $\|w(k)\| \leq d/2$. Define the augmented system $\Sigma_{\mathrm{a}}$ driven by $w$:
\begin{equation}\label{eq:augmented_lti}
\mathbf{A} := \begin{bmatrix} A & B\C \\ 0 & \A + \B_2 \C \end{bmatrix}, \quad
\mathbf{B}_2 := \begin{bmatrix} B \\ \B_2 \end{bmatrix}, \quad
\mathbf{C} := \begin{bmatrix} C & 0 \end{bmatrix}.
\end{equation}
If $\mathbf{A}$ is Schur stable, then
\begin{equation}\label{eq:EQ_lti_general}
E(Q) = \left\| \sum_{k=0}^{\infty} \mathrm{abs}\!\left(\mathbf{C}\,\mathbf{A}^k \mathbf{B}_2 \right) \right\| \cdot \frac{d}{2}.
\end{equation}
\end{lemma}

In (\ref{eq:augmented_lti}), $A$, $B$, $C$ are the plant matrices in (\ref{eq:lti_plant}), and the calligraphic $\A$, $\B_2$, $\C$ are the quantizer matrices in (\ref{eq:lti_quantizer}). The boldface matrices describe the augmented state $z := [(x - x_r)^\top,\, \xi^\top]^\top$ via $z(k+1) = \mathbf{A}\,z(k) + \mathbf{B}_2\,w(k)$, where $x$ and $x_r$ are the plant states driven through $Q$ and directly by $u$, respectively. The $(1,2)$ block $B\C$ of $\mathbf{A}$ is the coupling from the quantizer state $\xi$ to the plant state error.

For the optimal quantizer $Q^c$ in (\ref{eq:Qc}), $\A = A$, $\B_2 = B$, and $\C = -(CB)^{-1}CA$. Substituting these into (\ref{eq:augmented_lti}), the $(2,2)$ block $\A + \B_2\C$ of $\mathbf{A}$ becomes $F := A - B(CB)^{-1}CA$, and $\mathbf{C}\mathbf{A}^k\mathbf{B}_2 = CF^kB$. Since $CF = 0$, this vanishes for all $k \geq 1$, so that only the term $\mathbf{C}\mathbf{B}_2 = CB$ remains in (\ref{eq:EQ_lti_general}), in agreement with Lemma~\ref{lem:lti_optimal}. The matrix $F$ is the state matrix of the one-step-delayed inverse of $P$ when $CB$ is nonsingular \citep{silverman1969}. Moreover, $CF = 0$ yields the identity
\begin{equation}\label{eq:pencil_identity}
\begin{bmatrix} I & 0 \\ -C & zI \end{bmatrix}
\begin{bmatrix} zI - A & -B \\ C & 0 \end{bmatrix}
\begin{bmatrix} I & 0 \\ -(CB)^{-1}CA & I \end{bmatrix} 
= \begin{bmatrix} zI - F & -B \\ 0 & CB \end{bmatrix},
\end{equation}
whose determinant gives $z^m \det\left[\begin{smallmatrix} zI - A & -B \\ C & 0 \end{smallmatrix}\right] = \det(zI - F)\det(CB)$. The system pencil of $P$ thus has full normal rank, and the eigenvalues of $F$ are the $n - m$ invariant zeros of $P$ together with $m$ eigenvalues at the origin. Hence, $Q^c$ is stable if and only if all invariant zeros of $P$ lie in the open unit disk, i.e., $P$ is minimum phase. When $P$ is nonminimum phase, the linear programming-based approach of \citet{azuma2008tac} and the invariant set analysis approach of \citet{sawada2011acc} provide stable suboptimal quantizers with guaranteed performance bounds.
\subsection{Problem Formulation}\label{sec:problem}

\begin{figure}[t]
\centering
\includegraphics[width=0.8\linewidth]{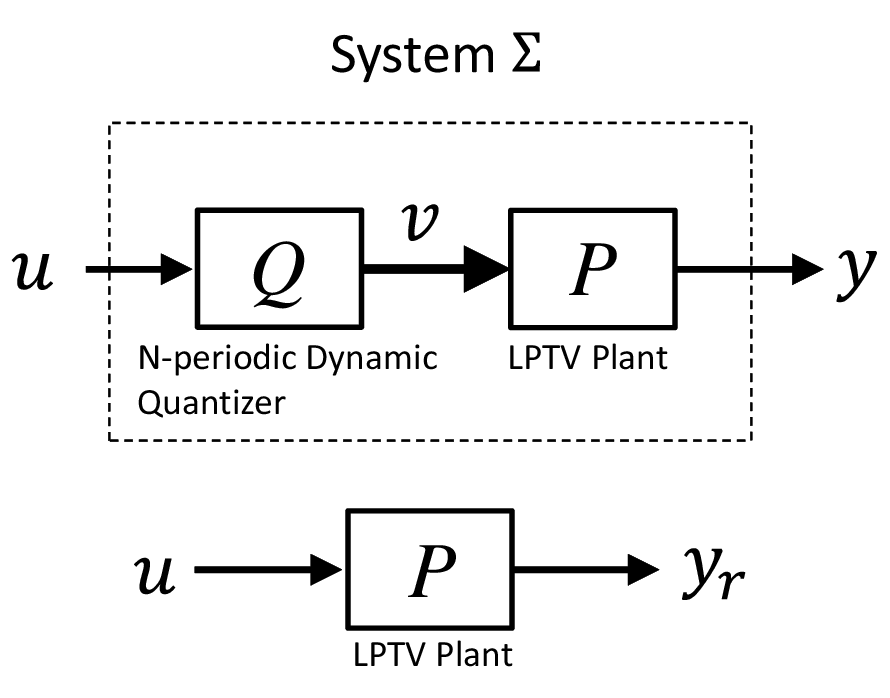}
\caption{Block diagram of the system $\Sigma$ (top) and the reference system (bottom). The dynamic quantizer $Q$ maps the continuous-valued input $u$ to the discrete-valued signal $v \in \mathbb{V}^m$, which drives the LPTV plant $P$. The reference output $y_r$ is obtained by applying $u$ directly to $P$.}
\label{fig:system}
\end{figure}

The system considered in this paper is illustrated in Fig.~\ref{fig:system}. The system $\Sigma$ consists of a dynamic quantizer $Q$ cascaded with an LPTV plant $P$. The quantizer $Q$ receives a continuous-valued input $u(k)$ and produces a discrete-valued output $v(k) \in \mathbb{V}^m$, which drives the plant $P$ to generate the output $y(k)$. As a reference, the output $y_r(k)$ denotes the response of the same plant $P$ driven directly by $u(k)$ without quantization. The design objective is to find $Q$ that minimizes the worst-case output difference between $y(k)$ and $y_r(k)$ over all admissible inputs.

The LPTV plant $P$ in (\ref{eq:lptv_state})--(\ref{eq:lptv_output}) is now driven by the discrete-valued control input $v \in \mathbb{V}^m$ rather than $u$:
\begin{equation}\label{eq:P_periodic}
P: \quad x(k+1) = A_k x(k) + B_k v(k), \quad y(k) = C_k x(k).
\end{equation}
The reference system is the same plant driven directly by $u$, i.e., $x_r(k+1) = A_k x_r(k) + B_k u(k)$ and $y_r(k) = C_k x_r(k)$, with the common initial state $x_r(0) = x(0) = x_0$. The dynamic quantizer for this periodic plant is allowed to be periodically time-varying with the same period $N$:
\begin{equation}\label{eq:Q_periodic}
Q: \begin{cases}
\xi(k+1) = \A_k \xi(k) + \B_{1,k} u(k) + \B_{2,k} v(k), \\
v(k) = q(\C_k \xi(k) + u(k)),
\end{cases}
\end{equation}
where $\xi \in \mathbb{R}^{n_\xi}$ is the quantizer state with $\xi(0) = 0$, $u \in \mathbb{R}^m$ is the continuous-valued input to the system $\Sigma$, and all coefficient matrices $\A_k$, $\B_{1,k}$, $\B_{2,k}$, $\C_k$ are $N$-periodic. The quantizer state dimension $n_\xi$ is a design parameter that may differ from the plant state dimension $n$; as shown in Theorem~\ref{thm:main}, the optimal quantizer has $n_\xi = n$.

Analogously to the LTI case (\ref{eq:E_def}), the maximum output difference between the quantized output $y(k)$ and the reference output $y_r(k)$ is defined as
\begin{equation}\label{eq:E_periodic}
E(Q) := \sup_{U \in l_\infty^m} \|Y_r(x_0, U) - Y(x_0, Q(U))\|,
\end{equation}
where $Y_r(x_0, U)$ denotes the output sequence of $P$ driven directly by $U$ (i.e., $y_r(k)$ in the lower part of Fig.~\ref{fig:system}), and $Y(x_0, Q(U))$ denotes the output sequence of $P$ driven through $Q$ (i.e., $y(k)$ in the upper part of Fig.~\ref{fig:system}). Since both outputs start from the same initial state $x_0$, their difference, and hence $E(Q)$, does not depend on $x_0$.

The following assumptions are imposed.
\begin{ass}\label{ass:dims}
The dimensions of $v$ and $y$ are the same ($m = p$).
\end{ass}

\begin{ass}\label{ass:quantizer}
For given quantization interval $d \in \mathbb{R}_+$, $\mathbb{V} := d\mathbb{Z}$, and the static quantizer $q$ is of the nearest-neighbor type satisfying $\mathrm{abs}(q(\mu) - \mu) \leq \frac{d}{2}\mathbf{1}_m$ for all $\mu \in \mathbb{R}^m$.
\end{ass}

\begin{ass}\label{ass:CB}
The matrix $C_{k+1} B_k$ is nonsingular for every $k \in \{0, 1, \ldots, N-1\}$.
\end{ass}

\begin{remark}\label{rem:ass3_meaning}
Assumption \ref{ass:CB} is the periodic counterpart of the condition that $CB$ is nonsingular in the LTI case. It requires that the system has relative degree one at every time instant: the input $v(k)$ directly affects the next output through $y(k+1) = C_{k+1}(A_k x(k) + B_k v(k))$ via the nonsingular gain $C_{k+1}B_k$. When $N = 1$, it reduces to $CB$ being nonsingular. The case where this assumption fails is discussed in Section \ref{sec:discussion}.
\end{remark}

\begin{prob}\label{prob:main}
For the system $\Sigma$ composed of the periodic plant $P$ in (\ref{eq:P_periodic}) and the periodic quantizer $Q$ in (\ref{eq:Q_periodic}), under Assumptions \ref{ass:dims}--\ref{ass:CB}, find $Q$ minimizing $E(Q)$ in (\ref{eq:E_periodic}) and determine the minimum value.
\end{prob}

In Section~\ref{sec:results}, the worst-case output error of a general periodic quantizer is first characterized. An optimal quantizer is then constructed by applying the LTI results of Section~\ref{sec:lti_review} to the cycled plant, and its optimality is verified directly in the original time domain.

\section{Main Results}\label{sec:results}

\subsection{Performance Analysis of Periodic Dynamic Quantizers}\label{sec:performance}

The LTI performance characterization (Lemma \ref{lem:lti_performance}) is first extended to the LPTV setting. Consider a quantizer $Q$ of the form (\ref{eq:Q_periodic}) with $\B_{1,k} = -\B_{2,k}$ for all $k$, and let $w(k) := v(k) - (\C_k\xi(k) + u(k))$ denote the static quantization error, which satisfies $\|w(k)\| \leq d/2$. Analogously to (\ref{eq:augmented_lti}), define the $N$-periodic matrices
\begin{align}\label{eq:augmented_ptv}
\mathbf{A}_k &:= \begin{bmatrix} A_k & B_k\C_k \\ 0 & \A_k + \B_{2,k}\C_k \end{bmatrix}, \quad
\mathbf{B}_{2,k} := \begin{bmatrix} B_k \\ \B_{2,k} \end{bmatrix}, \nonumber\\
\mathbf{C}_k &:= \begin{bmatrix} C_k & 0 \end{bmatrix},
\end{align}
with the same typographic convention as in the LTI case. The augmented state $z := [(x - x_r)^\top,\, \xi^\top]^\top$ obeys $z(k+1) = \mathbf{A}_k z(k) + \mathbf{B}_{2,k} w(k)$, where $x$ and $x_r$ are the plant states driven through $Q$ and directly by $u$, respectively, and the output error $e(k) := y_r(k) - y(k)$ is given by $e(k) = -\mathbf{C}_k z(k)$.

Since $\xi(0) = 0$ and $x(0) = x_r(0) = x_0$, it follows that $z(0) = 0$, and the output error at time $\tau$ is expressed as
\begin{equation}\label{eq:e_tau}
e(\tau) = -\sum_{l=1}^{\tau} \mathbf{C}_{\tau}\,\boldsymbol{\Psi}(\tau,\,\tau{-}l{+}1)\,\mathbf{B}_{2,\tau-l}\,w(\tau{-}l),
\end{equation}
where the state transition matrix is defined by
\begin{equation}\label{eq:aug_Psi}
\boldsymbol{\Psi}(j,i) := \begin{cases}
\mathbf{A}_{j-1}\mathbf{A}_{j-2}\cdots\mathbf{A}_{i}, & j > i, \\
I, & j = i.
\end{cases}
\end{equation}
The impulse response matrix characterizing the contribution of $w(\tau-l)$ to $e(\tau)$ at output time phase $i := \tau \bmod N$ is defined as
\begin{equation}\label{eq:Hli}
H_l(i) := \mathbf{C}_{i}\,\boldsymbol{\Psi}(i,\,i{-}l{+}1)\,\mathbf{B}_{2,i-l}, \quad l = 1, 2, \ldots,
\end{equation}
where all subscripts of $\mathbf{A}$ in $\boldsymbol{\Psi}(i,\,i{-}l{+}1)$ are reduced modulo $N$. Note that the input matrix $\mathbf{B}_{2,i-l}$ varies with $l$, reflecting the periodically changing disturbance input channel. In the LTI case, (\ref{eq:Hli}) reduces to $H_l = \mathbf{C}\mathbf{A}^{l-1}\mathbf{B}_2$.

\begin{theorem}\label{thm:performance}
For the system $\Sigma$ with a quantizer $Q$ of the form (\ref{eq:Q_periodic}) satisfying $\B_{1,k} = -\B_{2,k}$ for all $k$, suppose the monodromy matrix $\mathbf{A}_{N-1}\mathbf{A}_{N-2}\cdots\mathbf{A}_0$ of the augmented periodic system (\ref{eq:augmented_ptv}) is Schur stable. Then
\begin{equation}\label{eq:EQ_general}
E(Q) = \max_{i \in \{0,1,\ldots,N-1\}} \left\| \sum_{l=1}^{\infty} \mathrm{abs}\!\left(H_l(i)\right) \right\| \cdot \frac{d}{2}.
\end{equation}
\end{theorem}
\begin{proof}
By periodicity, (\ref{eq:e_tau}) reads $e(\tau) = -\sum_{l=1}^{\tau} H_l(i)\,w(\tau-l)$ with $i = \tau \bmod N$. The Schur stability of the monodromy matrix implies that $\|H_l(i)\|$ decays geometrically in $l$, so the series in (\ref{eq:EQ_general}) converges. Since $\|w(k)\| \leq d/2$, one has $\|e(\tau)\| \leq \|\sum_{l=1}^{\tau} \mathrm{abs}(H_l(i))\| \cdot d/2$ for every $\tau$, and hence $E(Q)$ does not exceed the right-hand side of (\ref{eq:EQ_general}). Conversely, for any given $\xi(k)$, the input $u(k)$ can be chosen such that $w(k)$ is arbitrarily close to any prescribed point of the box $[-d/2, d/2]^m$, because the argument of $q$ in (\ref{eq:Q_periodic}) ranges over $\mathbb{R}^m$ as $u(k)$ does. Fix a phase $i$ and a row index $j$ attaining the maximum absolute row sum of $\sum_{l=1}^{\infty} \mathrm{abs}(H_l(i))$. For $\tau$ with $\tau \bmod N = i$, choose $u(0), \ldots, u(\tau-1)$ successively so that each component of $w(\tau-l)$, $l = 1, \ldots, \tau$, approaches $d/2$ in magnitude with the sign opposite to that of the corresponding entry of the $j$-th row of $H_l(i)$. The $j$-th component of $e(\tau)$ then approaches $d/2$ times the $j$-th row sum of $\sum_{l=1}^{\tau} \mathrm{abs}(H_l(i))$, which converges to that of the infinite sum as $\tau \to \infty$. Hence $E(Q)$ is not smaller than the right-hand side of (\ref{eq:EQ_general}).
\end{proof}

Since $\mathbf{A}_k$ is block upper triangular, the hypothesis of Theorem~\ref{thm:performance} requires that both the plant monodromy matrix $A_{N-1}\cdots A_0$ and the monodromy matrix of $\A_k + \B_{2,k}\C_k$ be Schur stable. The right-hand side of (\ref{eq:EQ_general}) coincides with the LTI expression (\ref{eq:EQ_lti_general}) evaluated for the cycled augmented system, because the $i$-th block-row of $\check{\mathbf{C}}\,\check{\mathbf{A}}^{l-1}\check{\mathbf{B}}_2$ contains $H_l(i)$ as its only nonzero block, where $\check{\mathbf{A}}$, $\check{\mathbf{B}}_2$, $\check{\mathbf{C}}$ are constructed from $\mathbf{A}_k$, $\mathbf{B}_{2,k}$, $\mathbf{C}_k$ via (\ref{eq:checkA})--(\ref{eq:checkC}).

\subsection{Optimal Dynamic Quantizer for Periodic Systems}\label{sec:optimal}

The optimal quantizer is derived by applying Lemma~\ref{lem:lti_optimal} to the cycled plant $\check{P}$, i.e., the LTI system (\ref{eq:cycled_state})--(\ref{eq:cycled_output}) with the input $\check{u}$ replaced by $\check{v}$, and by extracting the periodic coefficient matrices from the resulting LTI quantizer.

\begin{theorem}\label{thm:main}
For the LPTV system $\Sigma$ under Assumptions \ref{ass:dims}--\ref{ass:CB}, every $N$-periodic dynamic quantizer $Q$ of the form (\ref{eq:Q_periodic}) satisfies
\begin{equation}\label{eq:lower_bound}
E(Q) \geq \max_{i \in \{0,1,\ldots,N-1\}} \left\|\mathrm{abs}\!\left(C_{i+1} B_i\right)\right\| \cdot \frac{d}{2}.
\end{equation}
Moreover, an optimal periodically time-varying dynamic quantizer is given by
\begin{equation}\label{eq:Qstar}
Q^*: \begin{cases}
\xi(k+1) = A_k \xi(k) + B_k v(k) - B_k u(k), \\
v(k) = q\!\left(-(C_{k+1} B_k)^{-1} C_{k+1} A_k \,\xi(k) + u(k)\right),
\end{cases}
\end{equation}
where $\xi \in \mathbb{R}^n$ and $\xi(0) = 0$, and the minimum value of $E(Q)$ is
\begin{equation}\label{eq:Estar}
E(Q^*) = \max_{i \in \{0,1,\ldots,N-1\}} \left\|\mathrm{abs}\!\left(C_{i+1} B_i\right)\right\| \cdot \frac{d}{2}.
\end{equation}
\end{theorem}

\begin{proof}
Steps~1 and~2 derive $Q^*$ from the LTI theory through cyclic reformulation; Steps~3 and~4 verify its optimality directly in the original time domain.

\textit{Step 1: Optimal LTI quantizer for the cycled plant.}
Direct computation from (\ref{eq:checkA})--(\ref{eq:checkC}) gives
\begin{eqnarray}
\check{C}\check{B} &=& \Pi \cdot \mathrm{diag}(C_1B_0, C_2B_1, \ldots, C_0B_{N-1}), \label{eq:checkCB}\\
\check{C}\check{A} &=& \Pi \cdot \mathrm{diag}(C_1A_0, C_2A_1, \ldots, C_0A_{N-1}), \label{eq:checkCA}
\end{eqnarray}
where $\Pi$ is the block permutation matrix mapping block-column $j$ to block-row $j+1$. Since $\Pi$ is nonsingular and each $C_{i+1}B_i$ is nonsingular by Assumption \ref{ass:CB}, $\check{C}\check{B}$ is nonsingular, and Lemma \ref{lem:lti_optimal} applied to $\check{P}$ gives the optimal LTI quantizer
\begin{equation}\label{eq:checkQstar}
\check{Q}^*: \begin{cases}
\check{\xi}(k+1) = \check{A}\check{\xi}(k) + \check{B}\check{v}(k) - \check{B}\check{u}(k), \\
\check{v}(k) = \check{q}\!\left(-(\check{C}\check{B})^{-1}\check{C}\check{A}\,\check{\xi}(k) + \check{u}(k)\right),
\end{cases}
\end{equation}
where $\check{q}: \mathbb{R}^{Nm} \to \mathbb{V}^{Nm}$ is the nearest-neighbor quantizer acting componentwise. From (\ref{eq:checkCB}) and (\ref{eq:checkCA}),
\begin{equation}\label{eq:checkCBinvCA}
(\check{C}\check{B})^{-1}\check{C}\check{A} = \mathrm{diag}\!\Big(
(C_1 B_0)^{-1} C_1 A_0, \
(C_2 B_1)^{-1} C_2 A_1, \
\ldots, \
(C_0 B_{N-1})^{-1} C_0 A_{N-1}
\Big),
\end{equation}
which is block-diagonal with $N$ blocks of size $m \times n$.

\textit{Step 2: Parameter extraction.}
The coefficient matrices of $\check{Q}^*$ retain the cyclic structure ($\check{A}$ and $\check{B}$ are block-circulant, and (\ref{eq:checkCBinvCA}) is block-diagonal), so that extracting their blocks gives (\ref{eq:Qstar}). Indeed, since $\check{q}$ acts componentwise and $q(0) = 0$, the output equation of (\ref{eq:checkQstar}) preserves the zero-block pattern of $\check{\xi}(k)$ and $\check{u}(k)$. Hence, as in Section~\ref{sec:cyclic}, the cycled input $\check{u}(k) = e_i \otimes u(k)$ with $\check{\xi}(0) = 0$ yields $\check{\xi}(k) = e_i \otimes \xi(k)$ and $\check{v}(k) = e_i \otimes v(k)$, $i = k \bmod N$, where $\xi$ and $v$ are the state and output of (\ref{eq:Qstar}). Therefore, $\check{Q}^*$ is the cycled representation of the $N$-periodic quantizer (\ref{eq:Qstar}), which has state dimension $n$; the quantizer $Q^*$ of the theorem is thus obtained.

\textit{Step 3: Performance of $Q^*$.}
Let $x$ and $x_r$ be the plant states driven through $Q^*$ and directly by $u$, respectively. By the state equation of (\ref{eq:Qstar}), $x - x_r$ and $\xi$ obey the same recursion with the same initial value; hence $\xi(k) = x(k) - x_r(k)$ and the output error is $e(k) := y_r(k) - y(k) = -C_k\xi(k)$. With $w(k)$ denoting the static quantization error, the output equation of (\ref{eq:Qstar}) gives $v(k) - u(k) = w(k) - (C_{k+1}B_k)^{-1}C_{k+1}A_k\,\xi(k)$. Then
\begin{align}\label{eq:e_onestep}
e(k+1) &= -C_{k+1}\left(A_k\xi(k) + B_k(v(k) - u(k))\right) \nonumber\\
&= -C_{k+1}B_k\,w(k).
\end{align}
Since $e(0) = 0$ and $\|w(k)\| \leq d/2$, it follows that $\|e(k+1)\| \leq \|\mathrm{abs}(C_{k+1}B_k)\| \cdot d/2$ for all $k$, so that $E(Q^*)$ does not exceed the right-hand side of (\ref{eq:Estar}).

\textit{Step 4: Lower bound.}
The argument of \citet[Theorem~1(i)]{azuma2008automatica} for the first time step is applied at each time phase. Consider an arbitrary $Q$ of the form (\ref{eq:Q_periodic}) and fix $i \in \{0, 1, \ldots, N-1\}$. Let $u(k) = 0$ for $k < i$. Since $\xi(0) = 0$ and $q(0) = 0$, it follows that $v(k) = 0$ and $\xi(k) = 0$ for $k < i$, and hence $\xi(i) = 0$ and $x(i) = x_r(i)$. Therefore, $v(i) = q(u(i))$ and
\begin{equation}\label{eq:e_general}
e(i+1) = -C_{i+1}B_i\left(q(u(i)) - u(i)\right).
\end{equation}
As $u(i)$ ranges over $\mathbb{R}^m$, the quantization error $q(u(i)) - u(i)$ comes arbitrarily close to every point of the box $[-d/2, d/2]^m$, so that $\sup_{u(i)} \|e(i+1)\| \geq \|\mathrm{abs}(C_{i+1}B_i)\| \cdot d/2$. Since $E(Q)$ is not smaller than this supremum for every $i$, (\ref{eq:lower_bound}) follows. Applying (\ref{eq:lower_bound}) to $Q^*$ and combining it with Step~3 proves (\ref{eq:Estar}); hence $Q^*$ is optimal.
\end{proof}

\begin{remark}[Interpretation]\label{rem:interpretation}
The quantizer $Q^*$ in (\ref{eq:Qstar}) has the following structure at each time $k$: (i) the state satisfies $\xi(k) = x(k) - x_r(k)$, i.e., it tracks the accumulated quantization effect in the plant state coordinates; (ii) the compensation term $-(C_{k+1}B_k)^{-1}C_{k+1}A_k\xi(k)$ cancels the effect of past errors on $y(k+1)$, so that, by (\ref{eq:e_onestep}), only the current quantization error reaches the next output, through the instantaneous gain $C_{k+1}B_k$. In terms of Theorem~\ref{thm:performance}, $H_1(i) = C_iB_{i-1}$ holds irrespective of the quantizer coefficients, whereas $Q^*$ achieves $H_l(i) = 0$ for all $l \geq 2$, in parallel with $\mathbf{C}\mathbf{A}^k\mathbf{B}_2 = 0$ for $k \geq 1$ in the LTI case of Section~\ref{sec:lti_review}. Only the term $H_1(i)$ thus remains in (\ref{eq:EQ_general}), in agreement with (\ref{eq:Estar}), and the performance is dictated by the time phase with the largest gain.
\end{remark}

\begin{remark}[Relation to the LTI result]\label{rem:lti_relation}
Theorem~\ref{thm:main} is the periodic counterpart of Lemma~\ref{lem:lti_optimal}, with the same three statements: a lower bound valid for every quantizer of the considered form, an optimal quantizer, and the minimum value. When $N = 1$, it reduces to Lemma~\ref{lem:lti_optimal}, and (\ref{eq:Qstar}) reduces to (\ref{eq:Qc}). The minimum (\ref{eq:Estar}) equals the optimum $\|\mathrm{abs}(\check{C}\check{B})\| \cdot d/2$ that Lemma~\ref{lem:lti_optimal} gives for the cycled plant $\check{P}$, since each block-row of $\check{C}\check{B}$ contains exactly one nonzero block. The lower bound, however, does not follow from Lemma~\ref{lem:lti_optimal} applied to $\check{P}$, whose worst case is taken over all inputs in $l_\infty^{Nm}$ rather than over cycled inputs only; it is therefore established directly in Step~4.
\end{remark}

\subsection{Stability of the Optimal Quantizer}\label{sec:stability}

The optimal quantizer $Q^*$ is not always stable. Substituting the output equation of (\ref{eq:Qstar}) into the state equation gives the error dynamics
\begin{equation}\label{eq:error_dynamics}
\xi(k+1) = F_k \,\xi(k) + B_k w(k),
\end{equation}
where
\begin{equation}\label{eq:Fk}
F_k := A_k - B_k (C_{k+1} B_k)^{-1} C_{k+1} A_k, \quad k = 0, 1, \ldots, N-1,
\end{equation}
and $w(k)$ is the static quantization error, which satisfies $\|w(k)\| \leq d/2$. The unforced part of (\ref{eq:error_dynamics}) is an $N$-periodic linear system with state matrices $F_k$, which is asymptotically stable if and only if the monodromy matrix
\begin{equation}\label{eq:monodromy}
\Phi := F_{N-1} F_{N-2} \cdots F_0
\end{equation}
is Schur stable \citep{cyc2}. Accordingly, $Q^*$ is said to be \emph{stable} if $\Phi$ is Schur stable, in which case the quantizer state $\xi(k)$ remains bounded because $w$ is bounded.

\begin{theorem}\label{thm:stability}
Under Assumptions \ref{ass:dims}--\ref{ass:CB}, the optimal quantizer $Q^*$ in (\ref{eq:Qstar}) is stable if and only if all invariant zeros of the cycled plant $\check{P}$ lie in the open unit disk.
\end{theorem}
\begin{proof}
Let $\check{F} := \check{A} - \check{B}(\check{C}\check{B})^{-1}\check{C}\check{A}$. Since $(\check{C}\check{B})^{-1}\check{C}\check{A}$ is block-diagonal with blocks $(C_{k+1}B_k)^{-1}C_{k+1}A_k$ by (\ref{eq:checkCBinvCA}), $\check{F}$ has the block-circulant structure of $\check{A}$ in (\ref{eq:checkA}) with each $A_k$ replaced by $F_k$. Consequently, $\check{F}^N = \mathrm{diag}(\Phi_0, \Phi_1, \ldots, \Phi_{N-1})$ with $\Phi_j := F_{j+N-1} \cdots F_{j+1} F_j$ and $\Phi_0 = \Phi$. The matrices $\Phi_j$ are cyclic permutations of the same product and therefore share the same nonzero eigenvalues, so that $\rho(\check{F})^N = \rho(\check{F}^N) = \rho(\Phi)$; that is, $\Phi$ is Schur stable if and only if $\check{F}$ is. On the other hand, $\check{C}\check{B}$ is nonsingular by (\ref{eq:checkCB}), and applying (\ref{eq:pencil_identity}) to the cycled plant $\check{P}$ shows that the eigenvalues of $\check{F}$ are the invariant zeros of $\check{P}$ together with $Nm$ eigenvalues at the origin. Hence $\check{F}$ is Schur stable if and only if all invariant zeros of $\check{P}$ lie in the open unit disk.
\end{proof}

\begin{remark}[Periodically minimum phase plants]\label{rem:min_phase}
A plant satisfying the zero condition of Theorem~\ref{thm:stability} is called periodically minimum phase \citep{okajima2026inv}. The invariant zeros of $\check{P}$ correspond to the periodic invariant zeros of $P$ \citep{cyc2} and coincide with the transmission zeros of $\check{P}$ when $\check{P}$ is reachable and observable.
\end{remark}

\begin{remark}[Connection to stable inversion]\label{rem:inverse_connection}
Inversion of LPTV plants was studied by \citet{perdon1992automatica}, and \citet{vanzundert2019ijc} constructed bounded noncausal inverses for nonminimum phase LPTV plants; an explicit periodic realization of the one-step-delayed inverse of $P$ was obtained via cyclic reformulation by \citet{okajima2026inv}. The matrix $F_k$ and the gain $-(C_{k+1}B_k)^{-1}C_{k+1}A_k$ of the compensation term in (\ref{eq:Qstar}) coincide with the state and output matrices of this realization, whose stability is characterized by the same zero condition as in Theorem~\ref{thm:stability}. The optimal quantizer thus performs a partial inversion of the plant at each time step: the compensation term applies the one-step-ahead inverse to the accumulated error state $\xi$ to cancel the output contribution of past errors.
\end{remark}

\subsection{Numerical Example}\label{sec:example}

Consider a 3-periodic ($N = 3$) SISO system with
\begin{gather}
A_0 = \begin{bmatrix} 0 & 1 \\ -0.5 & 0.8 \end{bmatrix}, \quad
A_1 = \begin{bmatrix} 0.3 & 0.5 \\ -0.4 & 0.6 \end{bmatrix}, \label{eq:ex1_A}\\
A_2 = \begin{bmatrix} 0.2 & 0.4 \\ -0.3 & 0.5 \end{bmatrix}, \label{eq:ex1_A2}
\end{gather}
\begin{equation}\label{eq:ex1_B}
B_0 = \begin{bmatrix} 0 \\ 1 \end{bmatrix}, \
B_1 = \begin{bmatrix} 0.5 \\ 1 \end{bmatrix}, \
B_2 = \begin{bmatrix} 0.5 \\ 1 \end{bmatrix},
\end{equation}
\begin{gather}
C_0 = \begin{bmatrix} 1 & 0 \end{bmatrix}, \quad C_1 = \begin{bmatrix} 0.5 & 1 \end{bmatrix}, \quad C_2 = \begin{bmatrix} 1 & 0.5 \end{bmatrix}, \label{eq:ex1_C}\\
d = 2. \nonumber
\end{gather}
The matrices $A_k$ and $B_k$ are those of the examples in \citet{okajima2026inv}. The plant monodromy matrix $A_2A_1A_0$ has spectral radius $\rho \approx 0.20$, so the plant is stable.

\textit{Verification of Assumption \ref{ass:CB}:}
$C_1 B_0 = 1$, $C_2 B_1 = 1$, and $C_0 B_2 = 0.5$, all nonzero.

\textit{Optimal quantizer performance:}
From (\ref{eq:Estar}), $E(Q^*) = \max\{|C_1 B_0|, |C_2 B_1|, |C_0 B_2|\} \cdot d/2 = \max\{1, 1, 0.5\} \cdot 1 = 1$.

\textit{Stability analysis:}
The state matrices (\ref{eq:Fk}) are
\begin{equation}\label{eq:F012}
F_0 = \begin{bmatrix} 0 & 1 \\ 0 & -0.5 \end{bmatrix}, \
F_1 = \begin{bmatrix} 0.25 & 0.1 \\ -0.5 & -0.2 \end{bmatrix}, \
F_2 = \begin{bmatrix} 0 & 0 \\ -0.7 & -0.3 \end{bmatrix},
\end{equation}
and the monodromy matrix is
\begin{equation}
\Phi = F_2 F_1 F_0 = \begin{bmatrix} 0 & 0 \\ 0 & -0.02 \end{bmatrix}.
\end{equation}
The eigenvalues of $\Phi$ are $\{0, -0.02\}$, so $Q^*$ is stable and, by Theorem~\ref{thm:stability}, the plant is periodically minimum phase. The error cancellation property $C_1 F_0 = 0$, $C_2 F_1 = 0$, and $C_0 F_2 = 0$ can also be verified directly.

\begin{figure}[t]
\centering
\includegraphics[width=0.9\textwidth]{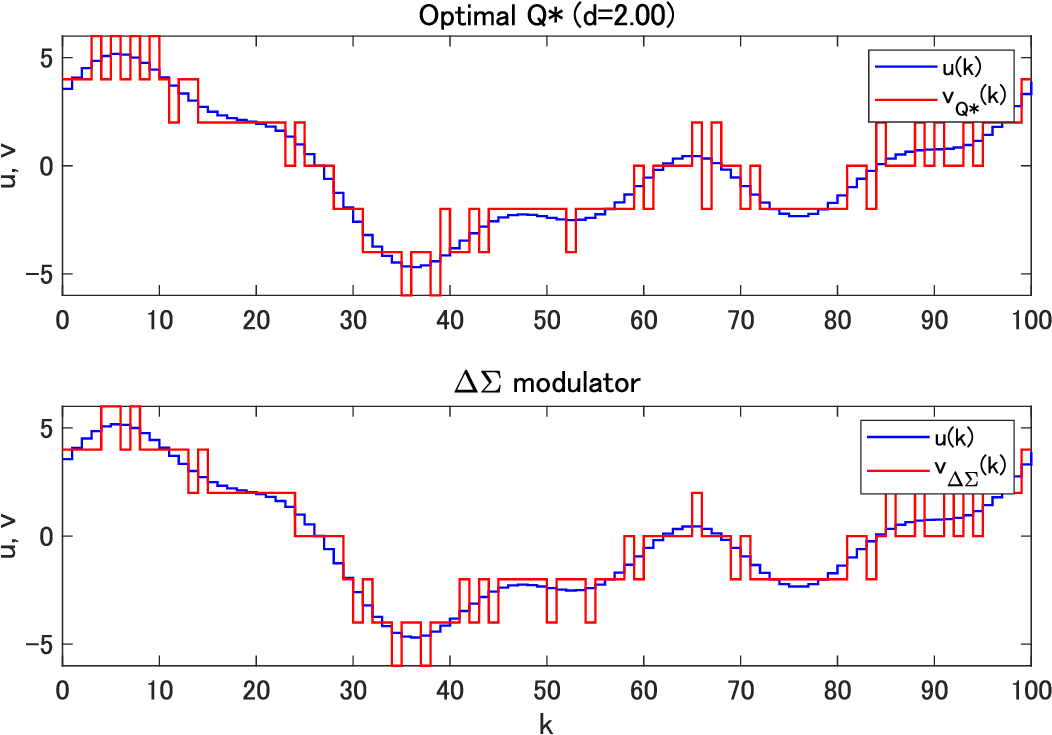}
\caption{Input $u(k)$ and quantized input $v(k)$ for the numerical example ($d = 2$). Top: optimal periodic $Q^*$. Bottom: $\Delta\Sigma$ modulator.}
\label{fig:inputs_d20}
\end{figure}

\textit{Simulation results:}
The proposed $Q^*$ is compared with two conventional quantizers: the first-order $\Delta\Sigma$ modulator \citep[see, e.g.,][]{quevedo2005tcsi} $v(k) = q(\xi(k) + u(k))$, $\xi(k+1) = \xi(k) + u(k) - v(k)$, a widely used dynamic quantizer that feeds back the accumulated quantization error without using any plant information, and the static quantizer $v(k) = q(u(k))$. The input signal is $u(k) = 3\cos(0.06(k+1)) + 2\sin(0.13(k+1)) + \sin(0.31(k+1))$ for $k = 0, 1, \ldots, 499$. Fig.~\ref{fig:inputs_d20} shows $u(k)$ and the quantized input $v(k)$ produced by $Q^*$ and by the $\Delta\Sigma$ modulator.

\begin{figure}[t]
\centering
\includegraphics[width=0.9\textwidth]{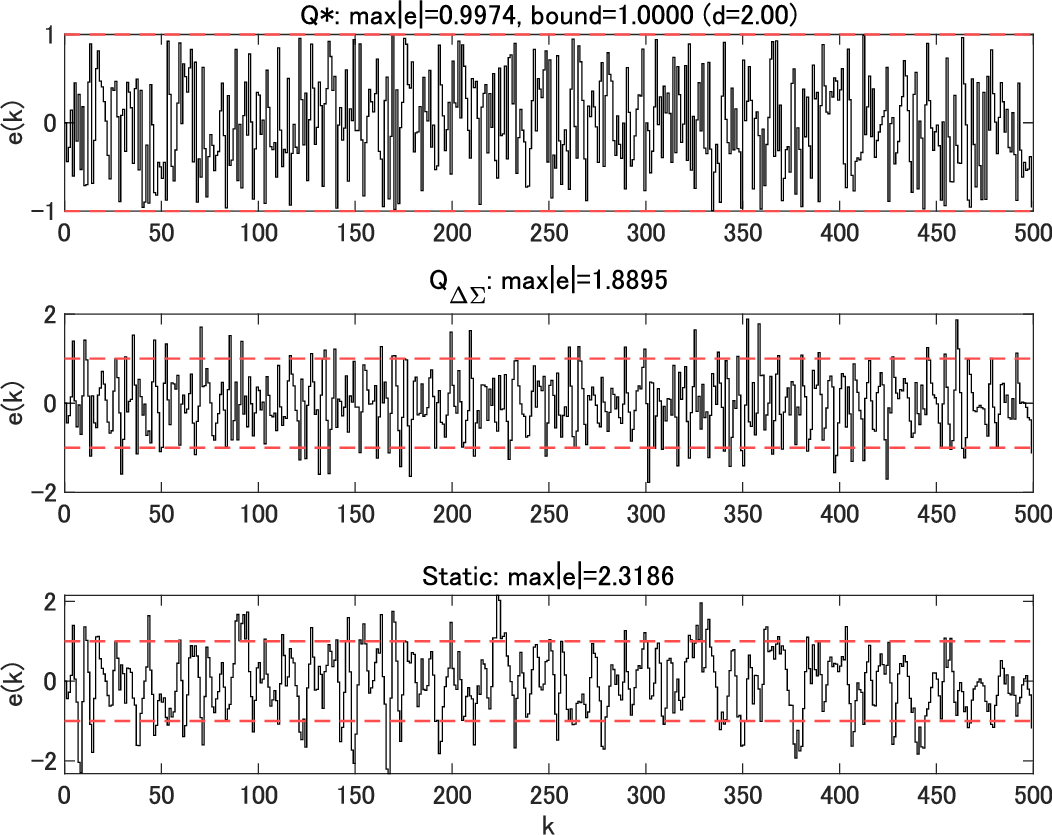}
\caption{Output errors for the numerical example ($d = 2$). Top: optimal periodic $Q^*$ (max$|e| = 0.997$, bound $= 1.0$). Middle: $\Delta\Sigma$ modulator (max$|e| = 1.889$). Bottom: static quantizer (max$|e| = 2.319$). Dashed lines: $\pm E(Q^*)$.}
\label{fig:errors_d20}
\end{figure}

Fig.~\ref{fig:errors_d20} shows the output error $e(k) = y_r(k) - y(k)$ for the three quantizers. The simulation is consistent with Theorem~\ref{thm:main}: $\max_k |e(k)| = 0.997$ for $Q^*$, which is below $E(Q^*) = 1$ (dashed lines). In contrast, the $\Delta\Sigma$ modulator yields $\max_k |e(k)| = 1.889$, and the static quantizer yields $\max_k |e(k)| = 2.319$.

\begin{figure}[t]
\centering
\includegraphics[width=0.9\textwidth]{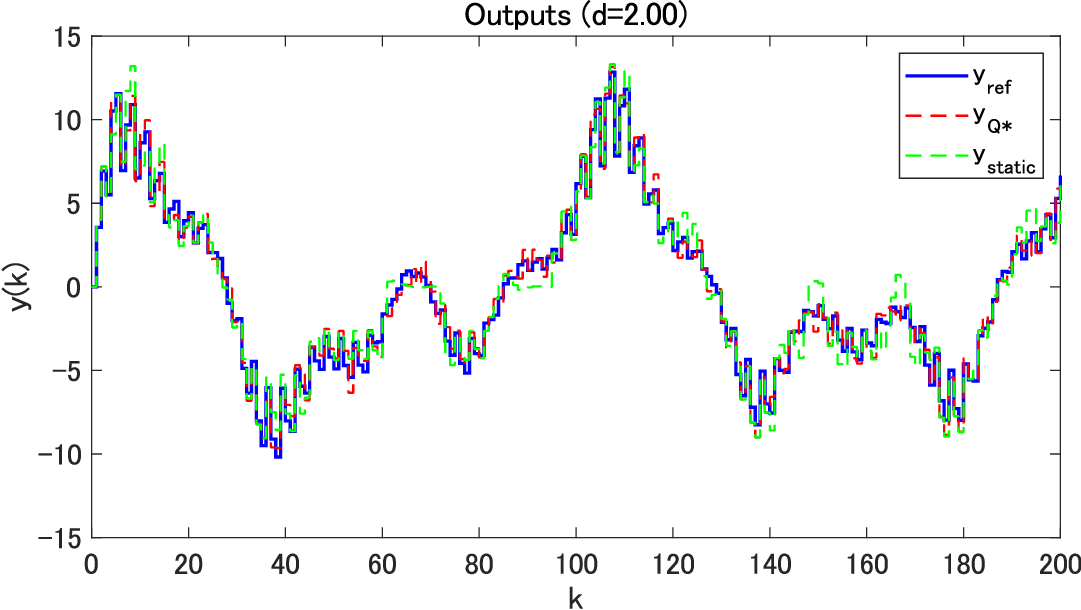}
\caption{Output comparison for the numerical example ($d = 2$). The output under $Q^*$ (red dashed) closely follows the reference $y_r$ (blue solid), while the static quantizer (green dashed) shows noticeable deviation.}
\label{fig:output_d20}
\end{figure}

Fig.~\ref{fig:output_d20} compares the plant output under the reference input $u(k)$ (no quantization, yielding $y_r(k)$), the optimal periodic quantizer $Q^*$, and the static quantizer. The output under $Q^*$ closely tracks the reference $y_r$, while the static quantizer exhibits visible deviation.

The maximum errors observed for a single input only bound $E(\cdot)$ from below. The worst-case errors of the two conventional quantizers are therefore evaluated by Theorem~\ref{thm:performance}, which applies to both: they are of the form (\ref{eq:Q_periodic}) with $\B_{1,k} = -\B_{2,k}$ (the $\Delta\Sigma$ modulator with $\A_k = 1$, $\B_{1,k} = 1$, $\B_{2,k} = -1$, $\C_k = 1$, and the static quantizer with all coefficients zero), and the corresponding augmented monodromy matrices are Schur stable because the plant is stable. Table~\ref{tab:comparison} lists the resulting values of $E(\cdot)$ together with the observed maximum errors, for $d = 2$ and $d = 0.5$; all values of $E(\cdot)$ scale linearly with $d$. In terms of the worst-case error, $Q^*$ achieves a reduction of about 65\% relative to both conventional quantizers, whereas the reduction observed in the simulation relative to the $\Delta\Sigma$ modulator (47--50\%) is specific to the input used. For this plant, the $\Delta\Sigma$ modulator, which uses no plant information, offers almost no worst-case improvement over the static quantizer; the improvement achieved by $Q^*$ stems from incorporating the plant dynamics into the quantizer.

\begin{table}[t]
\centering
\caption{Worst-case output error $E(\cdot)$ computed from (\ref{eq:EQ_general}) and maximum error observed in the simulation for the numerical example in Section~\ref{sec:example}.}
\label{tab:comparison}
\small
\begin{tabular}{llccc}
\toprule
& & $Q^*$ & $Q_{\Delta\Sigma}$ & Static \\
\midrule
$d = 2.0$ & $E(\cdot)$ & 1.000 & 2.870 & 2.885 \\
          & Observed   & 0.997 & 1.889 & 2.319 \\
$d = 0.5$ & $E(\cdot)$ & 0.250 & 0.717 & 0.721 \\
          & Observed   & 0.248 & 0.496 & 0.556 \\
\bottomrule
\end{tabular}
\end{table}

\section{Discussion}\label{sec:discussion}

The role of cyclic reformulation in the above results is twofold. First, it is constructive: the optimal LTI quantizer of the cycled plant retains the cyclic block structure, so that the periodic quantizer $Q^*$ with the original state dimension $n$ is read off block by block (Steps~1 and~2 of the proof of Theorem~\ref{thm:main}); no new optimization machinery is needed, and the optimality of $Q^*$ is then confirmed by an elementary per-step argument. Second, it links the stability of $Q^*$ to a time-invariant object, namely the invariant zeros of $\check{P}$ (Theorem~\ref{thm:stability}). Lifting \citep{chen1995sampled}, by contrast, stacks the inputs and outputs over one period and evolves on the $N$-step block time scale, so that the per-step causal structure required of a quantizer operating at the original sampling rate must be imposed separately.

The performance (\ref{eq:Estar}) shows that the achievable accuracy is dictated by the largest instantaneous gain $\max_i \|\mathrm{abs}(C_{i+1}B_i)\|$: a single time phase with large input-to-output coupling dominates the performance even if the gains at the other phases are small. Where the periodic structure is a design choice, it should therefore be chosen to reduce the largest of these gains. In addition, stability couples all $N$ time phases through the product $\Phi = F_{N-1}\cdots F_0$ (Section~\ref{sec:stability}). Implementing $Q^*$ requires only a few matrix--vector products of dimension $n$ and one static quantization per time step, together with the storage of the $N$ sets of coefficient matrices.

The present results require $C_{k+1}B_k$ to be nonsingular at every time phase (Assumption~\ref{ass:CB}) and, for stability, a periodically minimum phase plant (Theorem~\ref{thm:stability}). If the plant instead has uniform periodic relative degree $r \geq 2$, i.e., its first nonzero periodic Markov parameter is of order $r$ and nonsingular at every time phase \citep{okajima2026inv}, the first nonzero Markov parameter of the cycled plant is $\check{C}\check{A}^{r-1}\check{B}$, which is nonsingular and has the periodic Markov parameters as its only nonzero blocks; the closed-form design is then expected to carry over with $C_{k+1}B_k$ replaced by these parameters, which is left for future work. When the relative degree is not uniform, as in multirate input--output patterns where $C_{k+1}B_k = 0$ at some phases but not at others, $\check{C}\check{B}$ is singular without being zero, and no closed-form optimal quantizer is currently available. For nonminimum phase plants, the linear programming approach \citep{azuma2008tac} and the LMI-based invariant set approach \citep{sawada2011acc} can be applied to the LTI cycled plant $\check{P}$ to obtain stable suboptimal quantizers. However, the resulting quantizer has state dimension $Nn$ in general and is not guaranteed to retain the cyclic structure from which an $N$-periodic quantizer can be recovered, so that this structure has to be imposed in the synthesis.

\section{Conclusion}\label{sec:conclusion}

This paper has presented an optimal dynamic quantizer for square LPTV plants having relative degree one at every time phase. The quantizer $Q^*$ is obtained in closed form by applying the LTI theory to the cycled plant and reading off the periodic coefficient matrices block by block, and it is shown directly to attain the minimum worst-case output error, which is determined by the largest instantaneous gain $\max_i \|\mathrm{abs}(C_{i+1}B_i)\|$. The worst-case error of general periodic dynamic quantizers has also been characterized and used to quantify the advantage of $Q^*$ over conventional quantizers in a numerical example. The quantizer $Q^*$ embeds the one-step-delayed inverse of the plant, and it is stable if and only if the plant is periodically minimum phase, which is the condition for stable causal inversion of LPTV plants \citep{okajima2026inv}.

Future work includes the extension to feedback configurations, as done for LTI plants by \citet{minami2007cdc}, to plants with uniform periodic relative degree greater than one, and to plants with non-uniform relative degree, such as multirate systems in which Assumption~\ref{ass:CB} fails at some time phases.

\section*{Declaration of generative AI and AI-assisted technologies in the manuscript preparation process}
During the preparation of this work, the author used Claude (Anthropic) in order to refine English expressions and improve readability throughout the manuscript. After using this tool, the author reviewed and edited the content as needed and takes full responsibility for the content of the published article.

\section*{Funding}

This research did not receive any specific grant from funding agencies in the public, commercial, or not-for-profit sectors.


\end{document}